\documentclass[a4paper, 12pt, final]{elsarticle}

\usepackage{epsfig}
\usepackage{subfig}

\usepackage[hidelinks]{hyperref}
\usepackage{url}

\usepackage{amssymb}
\usepackage{amsthm}
\usepackage{amsmath}

\theoremstyle{plain}
\newtheorem{theorem}{Theorem}[section]

\newtheorem{corollary}[theorem]{Corollary}

\theoremstyle{definition}

\usepackage{booktabs} 
\usepackage{colortbl}
\usepackage[table]{xcolor}

\usepackage{tcolorbox}

\usepackage{lineno,hyperref}
\journal{Journal of Mathematical Analysis and Applications}

\begin{document}

\begin{frontmatter}

\title{When optimal is not the best: cost effectiveness analysis for HPV epidemic models}

\author[address1]{Fernando Salda\~na}
\ead{fernando.saldana@cimat.mx}

\author[address2]{Jos\'e Ariel Camacho-Guti\'errez}

\author[address3]{Ignacio Barradas}

\author[address4]{Andrei Korobeinikov}

\address[address1]{Instituto de Matem\'aticas, Campus Juriquilla, 76230, Universidad Nacional Aut\'onoma de M\'exico, Qu\'eretaro, Mexico}
\address[address3]{Centro de Investigaci\'{o}n en Matem\'{a}ticas, 36023 Guanajuato, Guanajuato, Mexico}
\address[address2]{Facultad de Ciencias, Universidad Aut\'onoma de Baja California, 22860 Baja California, Mexico}
\address[address4]{School of Mathematics and Information Science, Shaanxi Normal University, Xi\'an, China}

\begin{abstract}
This paper aims to evaluate the potential cost-effectiveness of healthcare interventions against human papillomavirus (HPV). For this, we consider a two-sex epidemic model for the transmission dynamics of HPV which includes screening, vaccination of adolescent boys and girls, and vaccination of sexually active adults. We first propose public health policies using constant control parameters and develop a cost-effectiveness analysis (CEA) to identify which intervention delivers the best effectiveness for the money invested. Secondly, we consider time-dependent control parameters and formulate an optimal control problem to obtain time-dependent versions of the interventions. As in the case of constant control parameters, we perform a CEA to investigate the cost-effectiveness of the time-dependent control interventions. Our findings suggest that females' vaccination, including adolescent girls and adult women, is the most cost-effective strategy. We also compare constant against the time-dependent healthcare interventions which are optimal in the sense that they minimize the objective functional of the optimal control problem. The results indicate that time-dependent controls are not always more cost-effective than constant controls. 
\end{abstract}


\end{frontmatter}

\section{Introduction}
%
%
%
Vaccination against the human papillomavirus (HPV) is the basic strategy for primary prevention of cervical cancer. The first vaccine to be approved to prevent HPV infection, Gardasil, has now been replaced by the 9-valent vaccine Gardasil-9. Gardasil-9 was originally approved for use in males and females aged 9 through 26 years. Nevertheless, in 2018, The U.S. Food and Drug Administration approved a supplemental application for Gardasil 9 expanding the approved use of the vaccine to include women and men aged 27 through 45 years.\par 
%
Since the introduction of HPV vaccines over a decade ago, HPV vaccination programs have been implemented in several countries. However, there are large discrepancies in coverage and targeted groups of vaccination strategies among countries \cite{soe2018should}. The vaccination program depends on country-specific factors, such as the economic and geographical constraints as well as the healthcare system organization. In several countries, prophylactic vaccination of pre-adolescent females has been introduced supported by modeled evaluations that have found this intervention to be cost-effective \cite{bogaards2015direct}. Vaccination of pre-adolescent males may also be cost-effective if females' coverage is below $50$\% \cite{canfell2012modeling}. Therefore, it is important to investigate under which conditions the inclusion of males and adult females into existing vaccination programs is cost-effective.\par 
Identifying interventions which have the potential to reduce substantially the disease burden at the best cost is of paramount importance for decision-makers faced with the choice of resource allocation. A number of mathematical models have been proposed to study the epidemiology of HPV and the development of vaccination strategies \cite{alsaleh2014, bogaards2011sex, guerrero2015cost, malik2013, smith2011predicted}. In this context, the optimal control theory has proven to be a very useful tool to identify optimal and potentially practical disease management strategies \cite{sharomi2017optimal}. In spite of this, to the best of authors' knowledge, only a very few studies have used the optimal control framework to investigate the cost-effectiveness of HPV vaccination programs \cite{brown2011role, malik2016optimal}.\par 
Here, we consider a two-sex epidemic model consisting of ordinary differential equations based on the susceptible-infected-vaccinated-susceptible (SIVS) compartmental structure. This model was considered in a previous work \cite{saldana2019optimal}, and the optimal, in the sense of minimizing the number of cervical cancer cases and related diseases in females was found there. In this study, we focus on the use of cost-effectiveness analysis to estimate the costs and health gains of alternative healthcare policies. Our ultimate goal is to prioritize the allocation of resources by identifying policies that have the potential to successfully control HPV infection in an efficient manner for the least resources.\par  
%
The content of the present paper can be outlined as follows. In the next section, we introduce the model and study its stability properties. In section \ref{sec:ICERconstant}, we propose realistic interventions based on the combination of the five constant control parameters incorporated into our model. Moreover, we define a suitable cost function and a cost-effectiveness ratio to evaluate the performance of these healthcare policies against HPV.
In section \ref{sec:OptimalControl}, we incorporate non-constant screening and vaccination rates into our model to analyze time-dependent healthcare policies via the optimal control theory. We also present a direct comparison between constant and time-dependent policies. Section \ref{sec:discussion} presents the discussion and concluding remarks.

\section{HPV model}\label{section:Model}
To study the transmission dynamics of the HPV in a heterosexual population, we consider a two-sex, compartmental model under the SIVS framework introduced in \cite{saldana2019optimal}. Only those HPV strains targeted by the nonavalent HPV vaccine are considered in this study. We classify the sexually active population of variable size $N(t)$ by the gender and the infection status; the subscripts $f$ and $m$ denote the female and male subpopulations, respectively.\par 
Since the majority of HPV cases are asymptomatic, we subdivide the female subpopulation of size $N_{f}(t)$ into four mutually exclusive compartments: the susceptibles, $S_{f}(t)$, the vaccinated, $V_{f}(t)$, and the infectious females unaware, $U_{f}(t)$, and aware, $I_{f}(t)$, of their infection, respectively. For the male subpopulation of size $N_{m}(t)$, we only consider three compartments: the susceptibles, $S_{m}(t)$, the vaccinated, $V_{m}(t)$, and the infectious males $I_{m}(t)$.\par
The basic assumptions of the model are as follows:
\begin{itemize}
\item[(i)] Individuals enter the sexually active population at a constant rate $\Lambda_{k}$, and leave the population by ceasing sexual activity at a per capita rate $\mu_{k}$ ($k=f,m$).
\item[(ii)] A fraction $w_{1}$ of females and a fraction $w_{2}$ of males (where $0\leq w_{1},w_{2}\leq 1$) are vaccinated before they enter the sexually active class and thus are recruited into their vaccinated compartment. Moreover, the susceptible sexually active females and males are vaccinated at per capita rates $u_{1}$ and $u_{2}$, respectively.
\item[(iii)] The vaccine reduces the force of infection by a factor $\epsilon$ with $0\leq \epsilon\leq 1$. Thus, the vaccine is absolutely effective when $\epsilon=0$ and absolutely ineffective when $\epsilon=1$. Therefore, $1-\epsilon$ is the vaccine efficacy. 
\item[(iv)] The exact duration of vaccine protection is unknown, but clinical trials \cite{huh2017} have shown sustained efficacy for at least $5$ years. We assume that vaccine-induced immunity wanes at a rate $\theta$; thus, for $\theta=0$ the protection is lifelong.
\item[(v)] Females have on average $c_{f}$ sexual contacts per unit of time. Thus, if $h$ is a transmission probability per contact, then the susceptible females are infected at a rate $c_{f}hI_{m}/N_{m}=\beta_{m}I_{m}/N_{m}$. Analogously, the susceptible males are infected by the unaware infected females at a rate $\beta_{f}U_{f}/N_{f}$ and by the aware infected females at a rate $\tilde{\beta}_{f}I_{f}/N_{f}$. We assume that $0<\tilde{\beta}_{f}<\beta_{f}$, since a female conscious of her infection may take precautions to reduce the probability of transmission. 
\item[(vi)] After infection, a fraction $p$ of the infected females may develop symptoms and become aware of their infection, entering the $I_{f}(t)$ class. The remaining infected females move to the unaware class $U_{f}(t)$. However, the screening allows the unaware infected females to detect their infection and change their status at a rate $\alpha$.
\item[(vii)] The infected individuals clear the infection naturally at a rate $\gamma_{k}$ ($k=f,m$). The existence and magnitude of the naturally acquired protection after HPV infection is still uncertain \cite{franceschi2014}. Therefore, we assume no permanent immunity after recovery. 
\end{itemize}\par
\begin{figure}[tbp]\centering
\includegraphics[width=\textwidth]{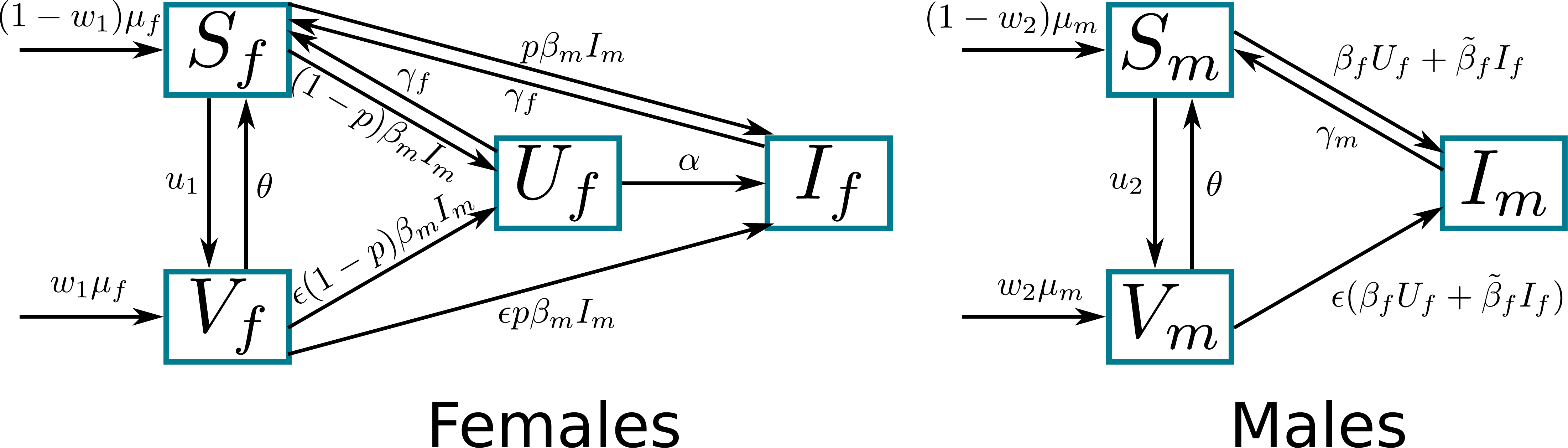}
\caption{Flow chart of model \eqref{ModelNormalized}.}
\label{fig:diagram}
\end{figure}

Moreover, since we are interested in the evaluation of several control interventions in a finite time interval $[0,T]$, the host population could be considered as constant. In mathematical terms, we assume $\Lambda_{k}=r_{k}N_{k}$ ($k=f,m$), where the per capita recruitment rate, $r_{k}$, is equal to the per capita rate of ceasing sexual activity, $\mu_{k}$. Please observe that $r_{k}N_{k}$ is the total number of individuals recruited into the sexually active population per unit of time. In addition, without loss of generality, we consider a rescaling of our model assuming that $N_{f}=1$ and $N_{m}=1$. Then, the state variables are expressed as fractions of the populations of each gender and we get the following system of differential equations:
\begin{equation}
\begin{aligned}
\dot{S}_{f}&=(1-w_{1})\mu_{f}-\beta_{m}S_{f}I_{m}-(u_{1}+\mu_{f})S_{f}+\gamma_{f}(U_{f}+I_{f})+\theta V_{f},\\
\dot{U}_{f}&=(S_{f}+\epsilon V_{f})(1-p)\beta_{m}I_{m}-(\gamma_{f}+\alpha+\mu_{f})U_{f},\\
\dot{I}_{f}&=(S_{f}+\epsilon V_{f})p\beta_{m}I_{m}+\alpha U_{f}-(\gamma_{f}+\mu_{f})I_{f},\\
\dot{V}_{f}&=w_{1}\mu_{f}+u_{1}S_{f}-\epsilon\beta_{m}V_{f}I_{m}-(\mu_{f}+\theta)V_{f},\\
\dot{S}_{m}&=(1-w_{2})\mu_{m}-(\beta_{f}U_{f}+\tilde{\beta}_{f}I_{f})S_{m}-(u_{2}+\mu_{m})S_{m}+\gamma_{m}I_{m}+\theta V_{m},\\
\dot{I}_{m}&=( \beta_{f}U_{f}+\tilde{\beta}_{f}I_{f})(S_{m}+\epsilon V_{m})-(\gamma_{m}+\mu_{m})I_{m},\\
\dot{V}_{m}&=w_{2}\mu_{m}-(\beta_{f}U_{f}+\tilde{\beta}_{f}I_{f})\epsilon V_{m}+u_{2}S_{m}-(\mu_{m}+\theta)V_{m}.
\label{ModelNormalized}
\end{aligned}
\end{equation} 
A flow chart depicting the model \eqref{ModelNormalized} is shown in Fig.~\ref{fig:diagram}. Here all the parameters are assumed to be nonnegative. The intervention measures, namely, the screening and vaccination rates will be called controls and denoted by vector $c=(w_{1},w_{2},u_{1},u_{2},\alpha)$. The biologically feasible region for system \eqref{ModelNormalized} is 
\begin{equation}
\Omega=\left\lbrace x\in \mathbb{R}_{+}^{7}: S_{f}+U_{f}+I_{f}+V_{f}= 1,\; S_{m}+I_{m}+V_{m}=1 \right\rbrace. 
\end{equation}\par 
Model \eqref{ModelNormalized} has a unique disease-free equilibrium state $E_{\circ}(c)$ which depends directly on the vector of controls $c$. The stability properties of $E_{\circ}(c)$ are determined by the effective reproduction number $\mathcal{R}_{e}(c)$.\par 
We remark that the definition of the basic reproduction number, $\mathcal{R}_{0}$, assumes a fully susceptible population and, hence, control measures such as mass vaccination that reduce the number of susceptible individuals in the population should, technically not reduce the value of $\mathcal{R}_{0}$. Therefore, alongside to $\mathcal{R}_{0}$ we introduce the effective reproduction number, denoted here as $\mathcal{R}_{e}$, which is defined as the actual average number of secondary cases per a primary case. $\mathcal{R}_{e}$ does not assume the complete susceptibility of the population; and, therefore, vaccination and other control measures could potentially reduce the value of $\mathcal{R}_{e}$. Consequently, in the presence of vaccination, the effective reproduction number can be a better metric for understanding the transmissibility of infectious diseases \cite{delamater2019}.\par 
%
Mathematically, both reproduction numbers can be computed via the next-generation operator introduced by Diekmann et al. \cite{diekmann1990definition}. Under this approach (see \cite{saldana2019optimal}), we obtain that the effective reproduction number is the geometric mean of the infection transfer from males to females $T_{m}^{f}(c)$ and the infection transfer from females to males $T_{f}^{m}(c)$, that is,
\begin{equation}
\mathcal{R}_{e}(c)=\sqrt{T_{m}^{f}(c)\cdot T_{f}^{m}(c)}.
\end{equation}\par
As a consequence of the van den Driessche theorem \cite{vanr0}, $\mathcal{R}_{e}(c)$ is a threshold value and we establish the following result regarding the local stability of the disease-free equilibrium.
\begin{corollary}\label{corollary}
The disease-free equilibrium $E_{\circ}(c)$ of system \eqref{ModelNormalized} is locally asymptotically stable for $\mathcal{R}_{e}(c)<1$ and unstable for $\mathcal{R}_{e}(c)>1$.
\end{corollary}\par
\begin{figure}[hbtp]
 \centering
  \subfloat[]{
    \includegraphics[width=0.5\textwidth]{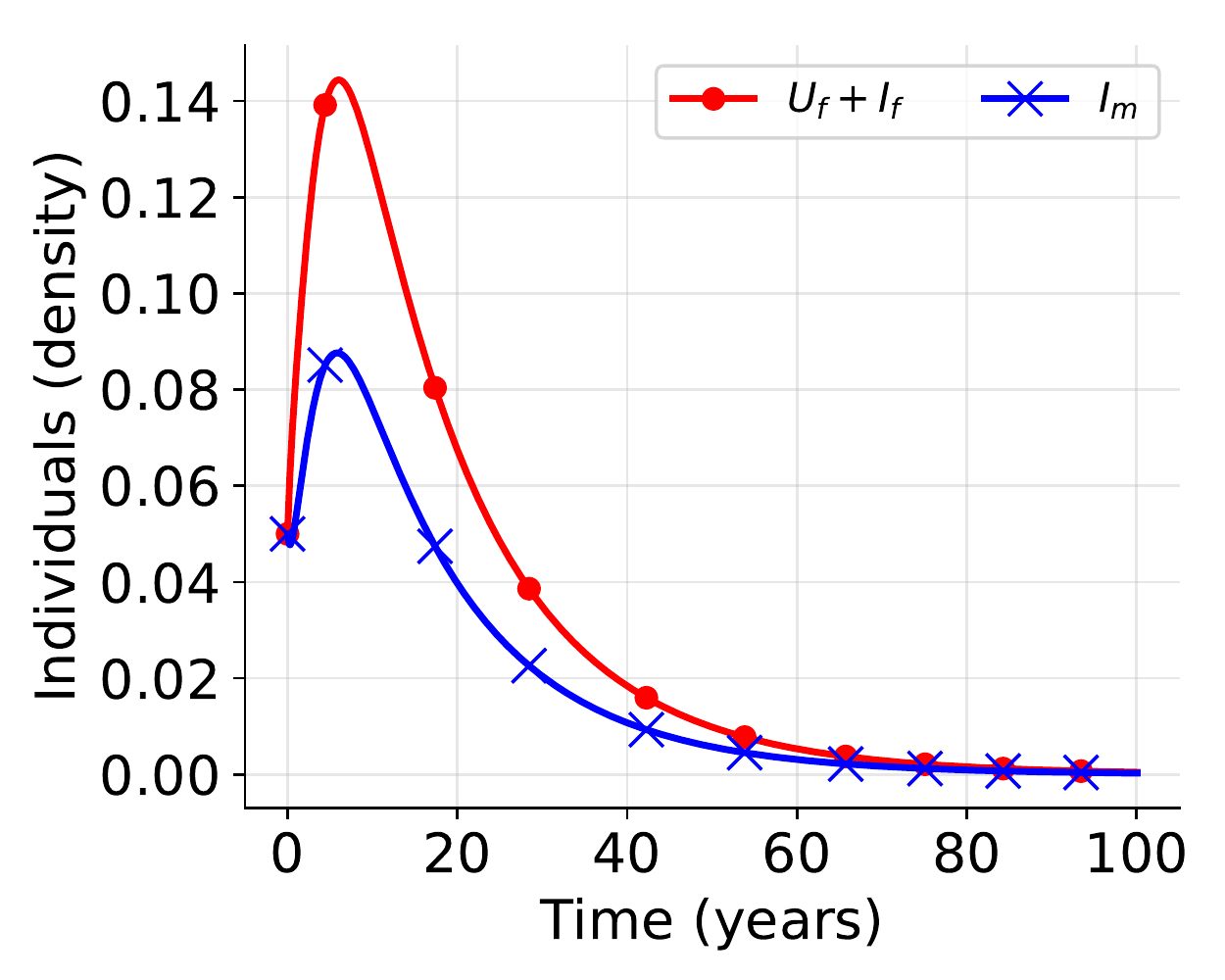}}
  \subfloat[]{
    \includegraphics[width=0.5\textwidth]{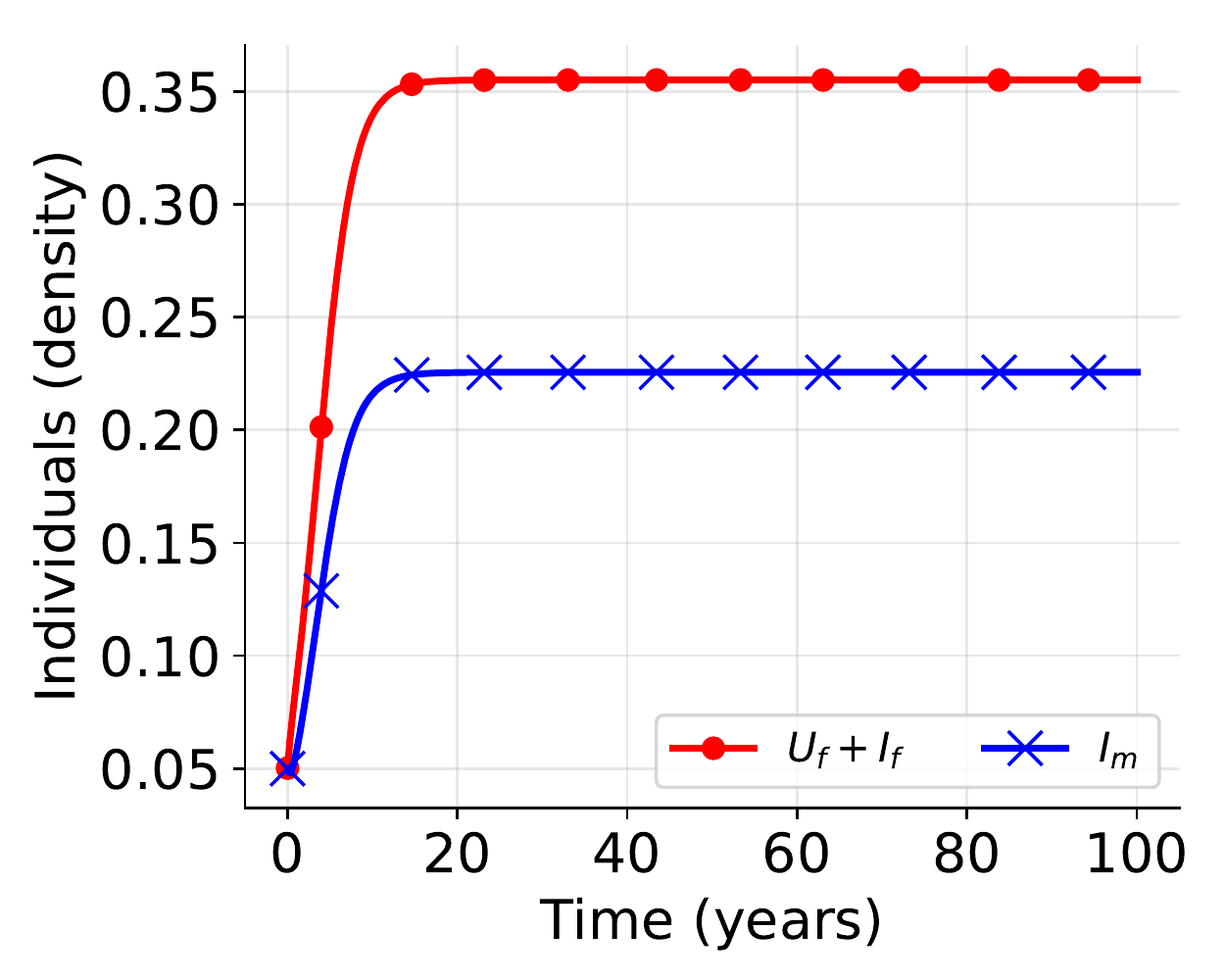}}
 \caption{Prevalence levels for model \eqref{ModelNormalized} with $\mathcal{R}_{e}(c)=0.9479<1$ (a) and $\mathcal{R}_{e}(c)=1.4151>1$ (b). In case (a) the values of the screening and vaccination rates are  $w_{1}=0.1$, $w_{2}=0.07$, $u_{1}=0.05$, $u_{2}=0.03$, and $\alpha=0.1$. In case (b), the screening and vaccination rates are taken equal to zero. The rest of the parameters are as described in Table \ref{TableParameters}. The initial conditions are: $S_{f}(0)=0.95$, $U_{f}(0)=0.03$, $I_{f}(0)=0.02$, $V_{f}(0)=0$, $S_{m}(0)=0.95$, $I_{m}(0)=0.05$, $V_{m}(0)=0$.}
 \label{fig:estability}
\end{figure}

 
Corollary \ref{corollary} implies that, if initial infection levels are sufficiently low, reducing and maintaining $\mathcal{R}_{e}(c)<1$ ensures disease's elimination. On the contrary, if the value of $\mathcal{R}_{e}(c)$ is higher than unity, the disease can persist in the population. Fig. \ref{fig:estability} illustrates these two possibilities. In Fig. \ref{fig:estability}(a), we choose screening and vaccination rates that guarantee $\mathcal{R}_{e}(c)<1$, so that the infection eventually decreases to zero. On the contrary, in Fig. \ref{fig:estability} (b), the screening and vaccination rates are taken equal to zero, thus giving $\mathcal{R}_{e}(0)>1$, so the system converges to an endemic equilibrium. Other parameters are fixed with their baseline values in Table \ref{TableParameters}.

\par


\begin{table}[hbtp]
\begin{center}
\begin{scriptsize}
\hspace*{-50pt}
\begin{tabular}{lcccc}
\toprule 
 \textbf{Parameter}  & \textbf{Range} & \textbf{Mean Value} & \textbf{Units} & \textbf{Source} \\
 \toprule 
Vaccine efficacy ($1-\epsilon$)& $[0.9, 1]$ & $0.95$ & \emph{adimensional} & \cite{brisson2016} \\[0.1cm]
 
Average duration of vaccine protection ($1/\theta$) & $[5, 50]$ & $20$ & \emph{year} & \cite{huh2017} \\[0.1cm]
 
 Male's transmission rate ($\beta_{m}$) & $[0.05, 5]$ & $2.0$ & \emph{year}$^{-1}$ & \cite{brisson2016}\\[0.1cm]
 
 Unaware females transmission rate ($\beta_{f}$) & $[0.05, 5]$ & $2.0$ & \emph{year}$^{-1}$ & \cite{brisson2016}\\[0.1cm]
 
 Aware females transmission rate ($\tilde{\beta}_{f}$) & $[0.025, 2.5]$ & $0.5$ & \emph{year}$^{-1}$ & Assumed \\[0.1cm]
 
 Female's infectious period ($1/\gamma_{f}$) & $[0.83, 2]$ & $1.3$  &  \emph{year} & \cite{munoz2004incidence} \\ [0.1cm]
 
 Male's infectious period ($1/\gamma_{m}$) & $[0.33, 1.2]$ & $0.6$  &  \emph{year} & \cite{anic2011} \\ [0.1cm]
 
 Female fraction that develop symptoms  ($p$) & $[0, 1]$ &  $0.4$ &\emph{adimensional} & Assumed \\[0.1cm]
 
Female's ceasing sexual activity rate ($\mu_{f}$) &  $[0.02, 1]$ & $1/20$ & \emph{year}$^{-1}$ & \cite{malik2013} \\[0.1cm]

 Male's ceasing sexual activity rate ($\mu_{m}$) & $[0.02, 1]$ & $1/25$ & \emph{year}$^{-1}$ & \cite{malik2013} \\ \bottomrule
 
\end{tabular}
\end{scriptsize}
\caption{Parameters of system \eqref{ModelNormalized}, sample units and source of estimation.} 
\label{TableParameters} 
\end{center}
\end{table}

\section{Cost-effectiveness analysis for constant controls}\label{sec:ICERconstant}


\begin{table}[t]

\centering
\begin{scriptsize}
\begin{tabular}{cl|c|c|c|c|c|}
\toprule 
ID &
Strategy description &
$w_1$ & $w_2$ & $u_1$ & $u_2$ & $\alpha$ \\ 
\toprule 
$S_1$ & All Controls 
& \checkmark & \checkmark 
& \checkmark & \checkmark 
& \checkmark \\ 
$S_2$ & Vaccination prior to sexual initiation
& \checkmark & \checkmark 
&  &  
&  \\ 
$S_3$ & Vaccination of sexually active individuals
&  & 
& \checkmark & \checkmark 
&  \\ 
$S_4$ & Females' vaccination 
& \checkmark & 
& \checkmark & 
&  \\ 
$S_5$ & Males' vaccination
&  & \checkmark 
&  & \checkmark 
&  \\ 
$S_6$ & Vaccination prior to sexual initiation and screening
& \checkmark & \checkmark 
&  &  
& \checkmark \\ 
$S_7$ & Vaccination of sexually active individuals and screening
&  & 
& \checkmark & \checkmark 
& \checkmark \\ 
$S_8$ &  Females' vaccination and screening
& \checkmark & 
& \checkmark & 
& \checkmark \\ \bottomrule
\end{tabular}

\caption{Control interventions analyzed in this work.}
\label{tab:strategies}
\end{scriptsize}
\end{table}
%
%
%
In the following, we carry out a cost-effectiveness analysis to investigate the most cost-effective control strategy against HPV transmission. We propose realistic strategies based on the combination of the five controls incorporated into our model. These strategies are shown in Table~\ref{tab:strategies}. For each strategy, the controls not marked with a checkmark are regarded as inactive.\par 
The main objective is to compare the health outcomes of the proposed interventions with respect to their application costs. To this end, we introduce the incremental cost-effectiveness ratio (ICER), which is usually defined as the additional cost, divided by its additional benefit compared with the next most expensive strategy. Mathematically speaking, for two strategies $S_1$ and $S_2$, the ICER is defined as
\begin{equation}
\label{ICER-two}
ICER(S_1, S_2) = \frac{C(S_2) - C(S_1)}{E(S_2) - E(S_1)},
\end{equation}
provided that $E(S_1) \neq E(S_2)$. Function $C(\cdot)$ measures the total costs (e.g. intervention cost, cost associated with illness, cost of treatment and monitoring, etc), while function $E(\cdot)$ measures the effectiveness of the intervention using some appropriate health outcome (such as lives saved, total number of infections averted, years of life gained, etc) \cite{okosun2011optimal}.\par 
We also define the average cost-effectiveness ratio (ACER) for a strategy $S$ as
\begin{equation}
\label{ICER-one}
ACER(S) = \frac{C(S)}{E(S)}.
\end{equation}\par 
The ICER could be interpreted as the additional cost incurred per additional health outcome, and the ACER as the net cost incurred per a unit of health outcome \cite{bang2014cost}.\par 
In this study, we measure the effectiveness of the intervention $S$ computing its cumulative level of infection averted using the following functional:
\begin{equation}
\label{CumuAverted}
E(S) = \int_0^T 
\left( U_f(t) - \tilde{U}_f(t) \right) +
\left( I_f(t) - \tilde{I}_f(t) \right) dt,
\end{equation}
where $U_f(t)$ are the infected unaware females for no control scenario at time $t$ and $\tilde{U}_f(t)$ are the infected unaware females under intervention $S$ (analogously for the infected aware females $I_f$). The functional \eqref{CumuAverted} focuses on the females because they are at a considerable higher risk of developing severe disease (e.g. cervical cancer) after an HPV infection.\par
To compute the costs, we propose the following functional:
\begin{tcolorbox}

\textbf{Cost functional}
\begin{eqnarray}
C(S)&=& \int_{0}^{T} L(S)+B_{1}\tilde{U}_{f}(t)+B_{2}\tilde{I}_{f}(t)dt,\label{CostC}
\end{eqnarray} 
where 
\begin{equation}
L(S)=A_{1}(w_{1}\mu_{f}+w_{2}\mu_{m})+A_{2}(u_{1}\tilde{S}_{f}(t)+u_{2}\tilde{S}_{m}(t))+A_{3}\alpha (\tilde{U}_{f}(t)+\tilde{S}_{f}(t)).
\end{equation}
\end{tcolorbox}
The parameters $A_{i}$ ($i=1,2,3$) are positive constants associated with the relative costs of vaccination and screening, and parameters $B_{j}$ ($j=1,2$) represent the medical and social costs associated with the illness. The terms $(w_{1}\mu_{f}+w_{2}\mu_{m})$ and $(u_{1}\tilde{S}_{f}(t)+u_{2}\tilde{S}_{m}(t))$ count the number of vaccinations given prior and after sexual initiation, respectively, and $\alpha (\tilde{U}_{f}(t)+\tilde{S}_{f}(t))$ counts the number of individuals that are screened.\par  
The costs attributable to the vaccination includes the price of the vaccine, the cost of delivery, as well as other administration associated costs. It is logical to expect that delivering the vaccine to school boys and girls is cheaper than vaccination of the sexually active individuals; hence, $A_{1}<A_{2}$. Since the costs of screening and follow up are uncertain, for simplicity, we assume that $A_{3}\approx A_{1}$; that is, the costs associated with screening and juveniles vaccination are of the same magnitude. Moreover, it is also considered that unaware infected females are at greater risk of developing HPV-induced cervical cancer than aware infected females; therefore, it can be expected that $B_{1}\geq B_{2}$. Under these considerations, as a base case, we assume $A_{1}=1$, $A_{2}=5$, $A_{3}=1$, $B_{1}=15$, $B_{2}=10$. Moreover, a time horizon of $T=100$ years is chosen because in this time frame the majority of the benefits and cost of vaccination can be recognized \cite{elbasha2010impact}. To perform a cost-effectiveness analysis, we use the following algorithm:

\textbf{\underline{Cost-Effectiveness Algorithm}}
\begin{enumerate}
\item The list of strategies with their corresponding costs and effectiveness is sorted from lowest to highest costs. Let $S_A$ and $S_B$ be the first and the second elements of the sorted list, respectively.

\item The $ACER(S_A)$ of the first element of the list $S_A$ is computed.

\item The $ICER(S_A, S_B)$ between the first two elements of the list $S_A, S_B$ is computed.

\item The following conditions are assessed:
\begin{itemize}
\item If $ICER(S_A, S_B) \leq 0$, then $S_A$ has higher effectiveness than $S_B$; hence, we keep $S_A$ and remove $S_B$.

\item If $ICER(S_A, S_B) \geq ACER(S_A)$, then $S_B$ has higher effectiveness than $S_A$, but in proportion is less cost-effective than $S_A$; hence, we keep $S_A$ and remove $S_B$.

\item If $0 < ICER(S_A, S_B) < ACER(S_A)$, then $S_B$ has higher effectiveness than $S_A$ and is in proportion more cost-effective than $S_A$; hence, we keep $S_B$ and remove $S_A$.
\end{itemize}

\item We return to STEP 2 until we find the most cost-effective strategy.
\end{enumerate}
To fairly compare the proposed strategies (see Table \ref{tab:strategies}), we set the values for the screening and vaccination rates in such a way that the value of the effective reproduction number $\mathcal{R}_{e}(c)$ is the same for all the strategies. Then, we compute the corresponding costs and the total infection averted over the time horizon for each intervention. These data allow us to compute both the ICER and the ACER using the algorithm proposed above and rank the interventions in order of increasing cost-effectiveness ratios. The results are summarized in Table \ref{tab:ICERresults}.\par    

%
\begin{table}[h!]
\centering
\begin{tabular}{lccc}
\toprule
Strategy &
$C(S_{i})$ &
$E(S_{i})$ &
Rank \\ 
\toprule
$S_1$ ($w_{1}=0.03$, $w_{2}=0.03$, $u_{1}=0.05$, $u_{2}=0.05$, $\alpha=0.1$) & \$ 70.33 & 31.77 & \# 7  \\ 
$S_2$ ($w_{1}=0.81$ $w_{2}=0.81$)                & \$ 47.86 & 31.04 & \# 2 \\ 
$S_3$ ($u_{1}=0.068$, $u_{2}=0.05$)              & \$ 69.07 & 31.71 & \# 6 \\ 
$S_4$ ($w_{1}=0.3$, $u_{1}=0.127$)               & \$ 49.24 & 32.43 & \# 1 \\ 
$S_5$ ($w_{2}=0.3$, $u_{2}=0.119$)               & \$ 55.07 & 31.86 & \# 3 \\ 
$S_6$ ($w_{1}=0.66$ $w_{2}=0.6$, $\alpha=0.4$)   & \$ 59.50 & 31.99 & \# 5 \\ 
$S_7$ ($u_{1}=0.046$, $u_{2}=0.05$, $\alpha=0.2$) &\$ 73.30 & 32.01 & \# 8 \\ 
$S_8$ ($w_{1}=0.15$, $u_{1}=0.1$, $\alpha=0.3$)   & \$ 58.03 & 32.65 & \# 4 \\ \bottomrule
\end{tabular}
\caption{Control strategies with fixed constant controls together with their costs (using \eqref{CostC}) and cumulative level of infection averted. The ranking of the strategies is according to the cost-effectiveness algorithm. For all the strategies, the value of the effective reproduction number is $\mathcal{R}_{e}(c)=0.9$.}
\label{tab:ICERresults}
\end{table}
The results of the cost-effectiveness analysis suggest that the strategy $S_{4}$ (female's vaccination) is the most cost-effective intervention and $S_{2}$ (vaccination before sexual initiation) is the strategy with the second-best performance. These results coincide with several studies that have analyzed HPV transmission at the population level finding that vaccination of pre-adolescent girls is both highly effective and highly cost-effective to reduce the disease burden caused by HPV, see Seto et al. \cite{seto2012cost} and the references therein. However, regardless of the independence of the cost function \eqref{CostC} from infected males, the third most cost-effective strategy is $S_{5}$ (male's vaccination). This result is somehow unexpected; however, it could be explained by the fact that for a heterosexual population, the eradication of the infection can be achieved by vaccinating a considerable proportion of a single-sex. On the other hand, the strategy with the worst performance is $S_{7}$ (vaccination of sexually active individuals with female's screening), thus, even assuming a medium duration of protection for the vaccine ($20$ years), it is plausible to vaccinate individuals before sexual debut avoiding the potential of an HPV infection.\par
\begin{figure}[h!]
\centering
\includegraphics[width=0.5\textwidth]{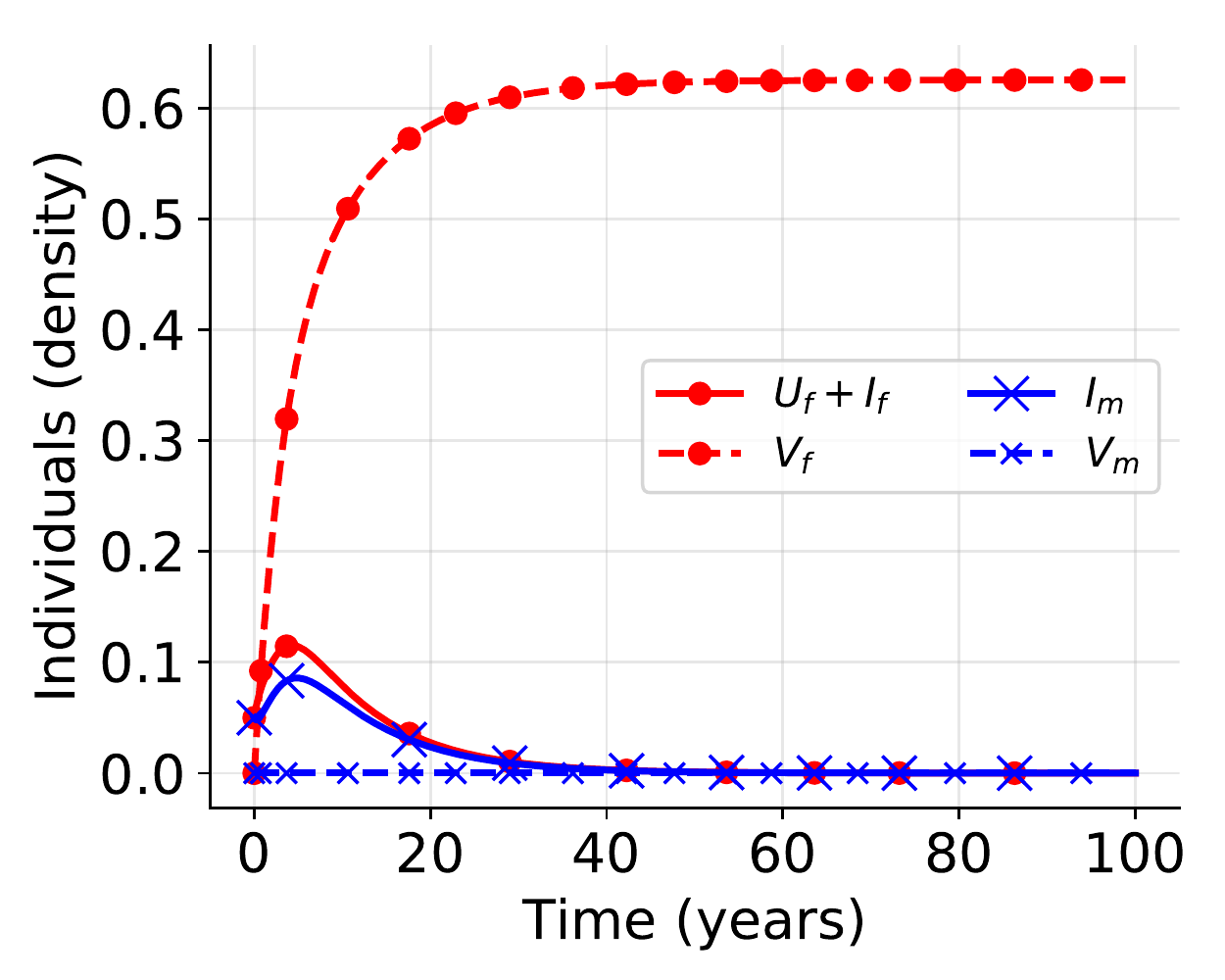}
\caption{Dynamics of model \eqref{ModelNormalized} state variables under intervention $S_{4}$ with fixed constant controls $w_1=0.3$ and $u_1=0.127$}
\label{fig:best_constant}
\end{figure}

The simulations of the HPV model showing the effects of strategy $S_{4}$, which is the most cost-effective control strategy according to the results in Table \ref{tab:ICERresults}, are illustrated in Fig. \ref{fig:best_constant}. The simulations show that the number of infected individuals decreases to zero after an approximate time of $40$ years. Furthermore, the fraction of vaccinated females increases until it reaches a value of $0.62$, whereas the fraction of vaccinated males is maintained at zero. 

\section{Cost-effectiveness analysis for optimal controls}\label{sec:OptimalControl}
In the previous sections, we analyzed the qualitative properties of the HPV model \eqref{ModelNormalized} and performed a cost-effectiveness analysis for constant control policies. Constant control analysis is helpful to understand the mean behavior of the model. However, such a constant control scheme disregards the changing dynamics of the infection and thereby can be a fragile strategy.\par
For monetary and time resources, as well as disease elimination, it is essential to find the right time and the right amount of control administration. Therefore, it is valuable to incorporate time-dependent screening and vaccination into our model and explore the most cost-effective strategy via the optimal control theory.\par 
Including time-dependent controls $c(t) = (w_1(t), w_2(t), u_1(t), u_2(t), \alpha(t))$ and substituting $S_{f}=1-U_{f}-I_{f}-V_{f}$, $S_{m}=1-I_{m}-V_{m}$, we obtain the following five-dimensional control model for HPV dynamics:
\begin{equation}
\begin{aligned}
\dot{U}_{f}&=\left((1-U_{f}-I_{f}-V_{f})+\epsilon V_{f}\right)(1-p)\beta_{m}I_{m}-(\gamma_{f}+\alpha(t)+\mu_{f})U_{f},\\
\dot{I}_{f}&=\left((1-U_{f}-I_{f}-V_{f})+\epsilon V_{f}\right)p\beta_{m}I_{m}+\alpha(t)U_{f}-(\gamma_{f}+\mu_{f})I_{f},\\
\dot{V}_{f}&=w_{1}(t)\mu_{f}+u_{1}(t)(1-U_{f}-I_{f}-V_{f})-\epsilon\beta_{m}V_{f}I_{m}-(\mu_{f}+\theta)V_{f},\\
\dot{I}_{m}&=(\beta_{f}U_{f}+\tilde{\beta}_{f}I_{f})((1-I_{m}-V_{m})+\epsilon V_{m})-(\gamma_{m}+\mu_{m})I_{m},\\
\dot{V}_{m}&=w_{2}(t)\mu_{m}-(\beta_{f}U_{f}+\tilde{\beta}_{f}I_{f})\epsilon V_{m}+u_{2}(t)(1-I_{m}-V_{m})-(\mu_{m}+\theta)V_{m}.
\label{ControlModel}
\end{aligned}
\end{equation} 
Here, the controls are subject to constraints
\begin{equation}
0\leq w_{1}(t), w_{2}(t), u_{1}(t), u_{2}(t), \alpha(t)\leq 1.\label{ConstrainsControls}
\end{equation}
The control model \eqref{ControlModel} is defined on a finite time interval $[0,T]$ and the set of admissible controls, denoted $D(T)$, is defined as the set of all possible Lebesgue measurable functions which for almost all $t\in[0,T]$ satisfy constraints \eqref{ConstrainsControls}.\par 
For an optimal control problem, it is necessary to define a quantitative criterion to evaluate the performance of the admissible controls. In mathematical terms, this quantitative criterion leads to the definition of an objective functional, and an optimal control is one that minimizes (or maximizes) this functional.
For the control model \eqref{ControlModel} and the set of admissible controls $D(T)$, we want to minimize the overall impact of the infection and the level of efforts that would be needed to control the infection over $T$ years. Thus, we consider system \eqref{ControlModel} together with the following objective functional
\begin{tcolorbox}

\textbf{Objective functional}
\begin{equation}
J(c)=\int_{0}^{T} B_{1}U_{f}+B_{2}I_{f}+\dfrac{1}{2}\left( A_{1}(w_{1}^{2}+w_{2}^{2})+A_{2}(u_{1}^{2}+u_{2}^{2})+A_{3}\alpha^{2}\right)dt.\label{ObjectiveFunctional}
\end{equation}
\end{tcolorbox}
It is generally assumed that the quadratic terms penalize high levels of control administration (to avoid costly interventions) because "increasing the availability of vaccines and other resources often leads to a waste" \cite{sharomi2017optimal}. 
Nevertheless, it is important to mention that the sum of weighted squares of the controls is probably the most common form of the objective functional in the literature due to its mathematical convenience. In particular, for functionals of these type, it is possible, by virtue of the Pontryagin maximum principle, to obtain the optimal controls as explicit functions of the state and adjoint variables.\par
For the sake of simplicity and comparison, in this paper the objective functional \eqref{ObjectiveFunctional} considers the sum of weighted squares of the controls and also uses the same weight parameters of the cost function \eqref{CostC}. 

\begin{tcolorbox}

\textbf{Optimal control problem}

The general optimal control problem is to find optimal vaccination and screening rates $c^{*}=(w_{1}^{*}(t), w_{2}^{*}(t), u_{1}^{*}(t), u_{2}^{*}(t), \alpha^{*}(t))$ such that
\begin{equation}
J(c^{*})=\min_{c\in D(T)} J(c)\label{OptimalControlProblem}
\end{equation}
subject to the dynamics of the HPV control model \eqref{ControlModel}. 
\end{tcolorbox}

Theorem 4.1 in \citep[Chapter~III]{Fleming1975} ensures the existence of an optimal control and the corresponding solution for this problem. Proofs of such statements can be found in \cite{camacho2018, saldana2019optimal}.\par 
Here, we must clarify the following: (i) Although optimal controls $c^{*}$ are expected to be more cost-effective than constant controls, this is not known \emph{a priori} because $c^{*}$ minimize the objective functional \eqref{ObjectiveFunctional} and the ICER measures the cost using \eqref{CostC}; (ii) The cost function \eqref{CostC} can also be used as an objective functional for the optimal control; however, generally speaking, an optimal control problem could have multiple goals and therefore, the objective functional may include other factors besides costs. 

\subsection{Characterization of the optimal controls}
Here, we obtain the optimality system that corresponds to complement the control model \eqref{ControlModel} with a dual system for adjoint variables. Then, we can achieve the characterization of the optimal controls in terms of the state and adjoint variables.\par
\begin{figure}[hbtp]
\centering
\subfloat[Strategy $S_{1}^{*}(t)$]{
\includegraphics[width=0.35\textwidth]{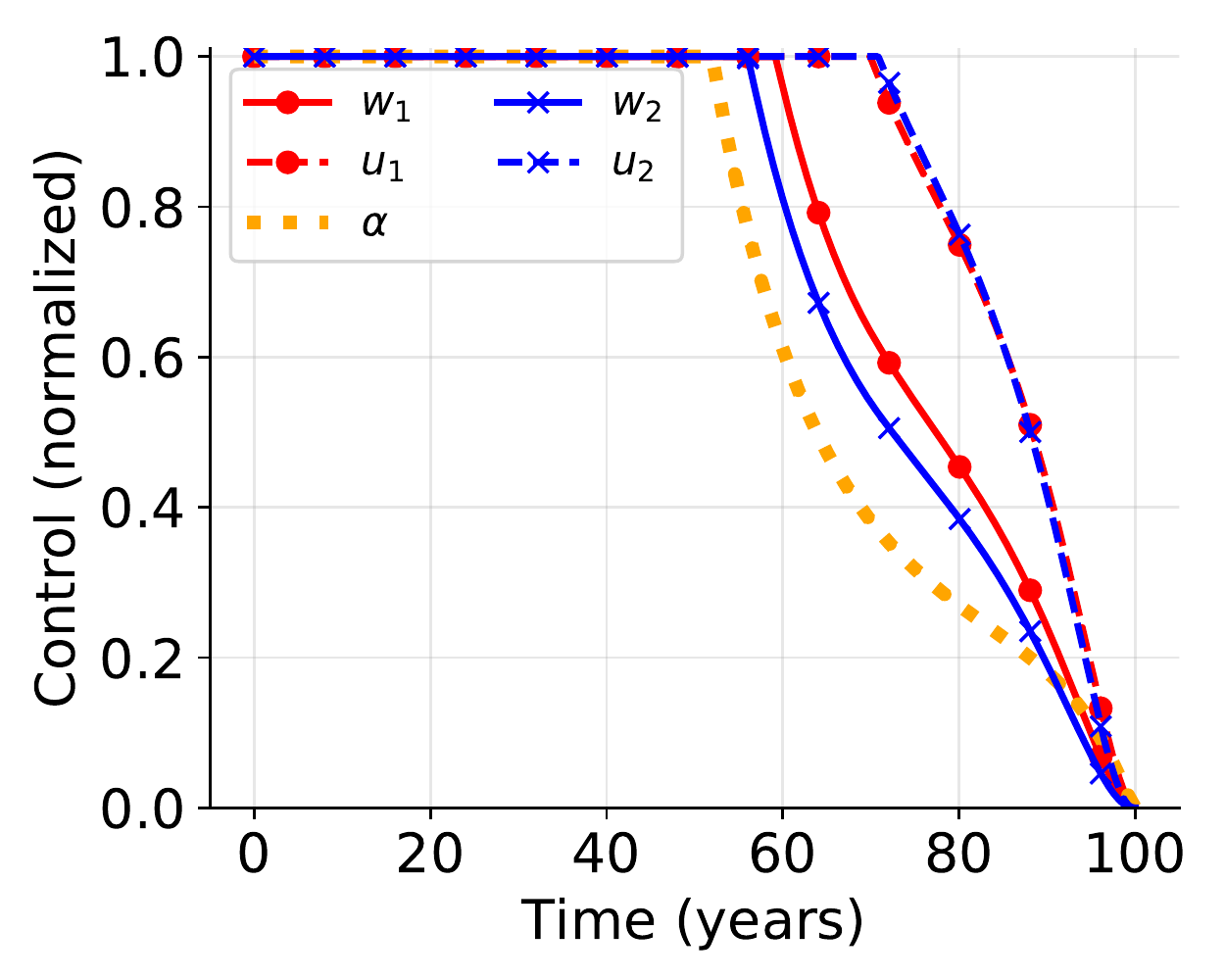}}
\subfloat[Strategy $S_{2}^{*}(t)$]{
\includegraphics[width=0.35\textwidth]{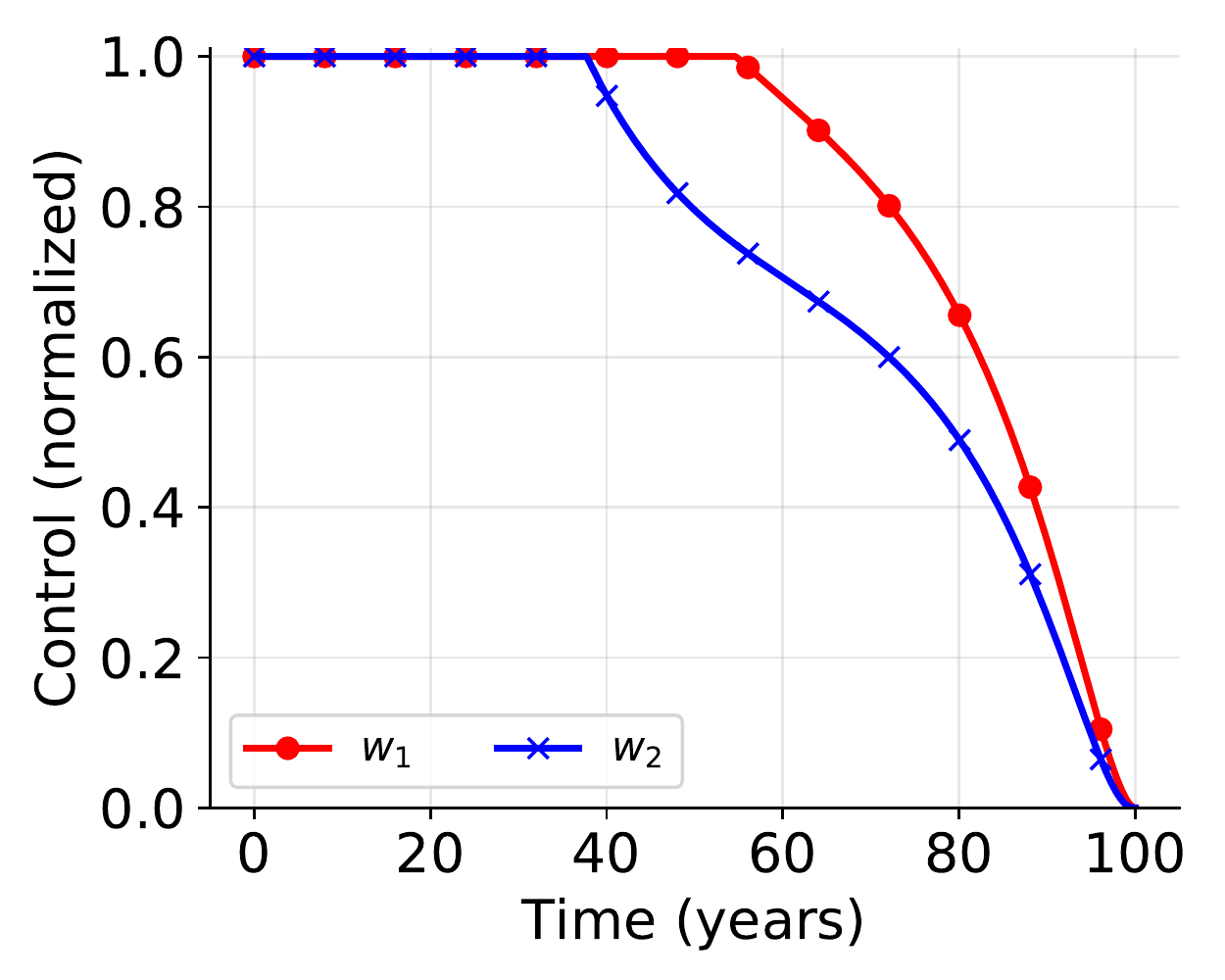}}

\subfloat[Strategy $S_{3}^{*}(t)$]{
\includegraphics[width=0.35\textwidth]{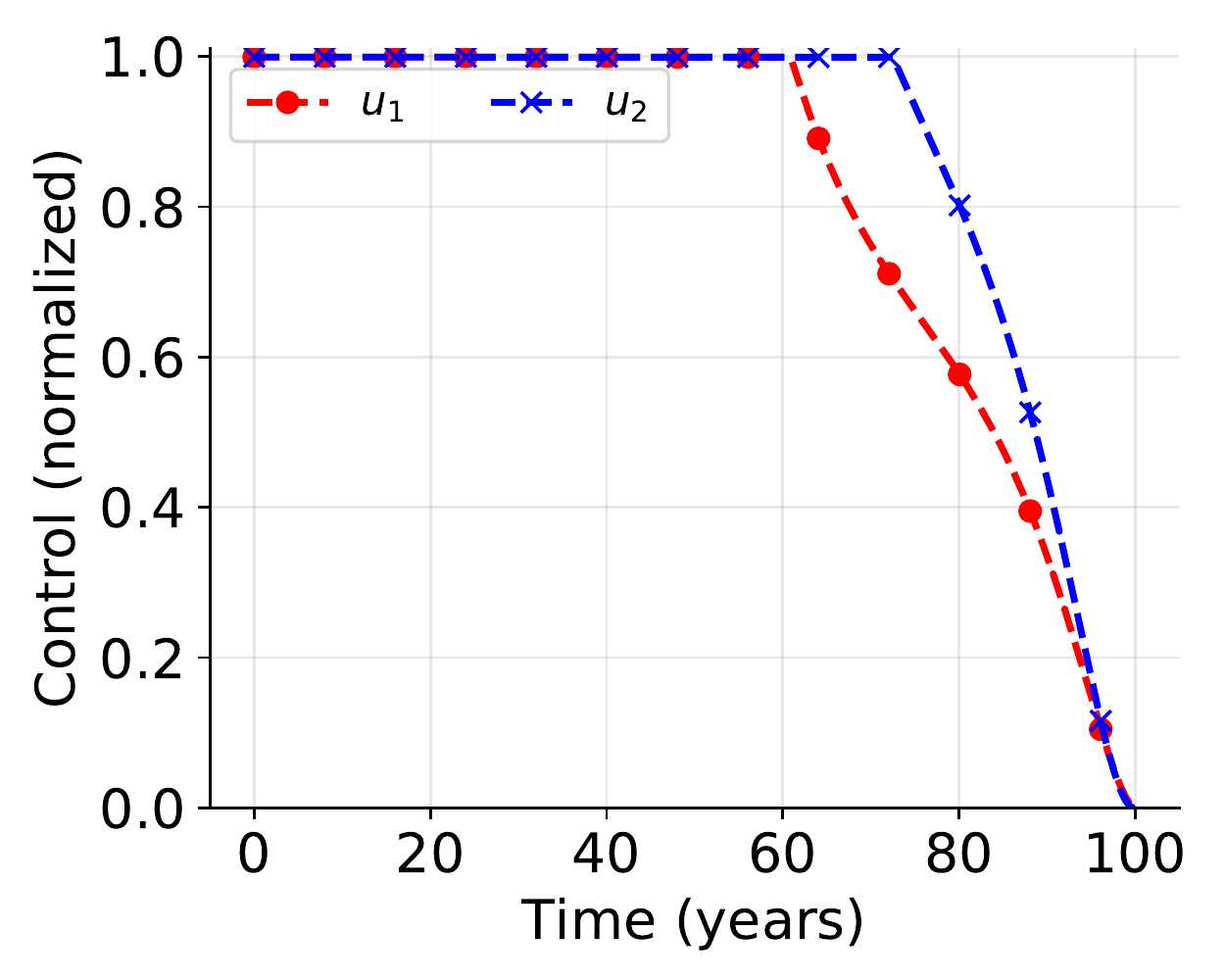}}
\subfloat[Strategy $S_{4}^{*}(t)$]{
\includegraphics[width=0.35\textwidth]{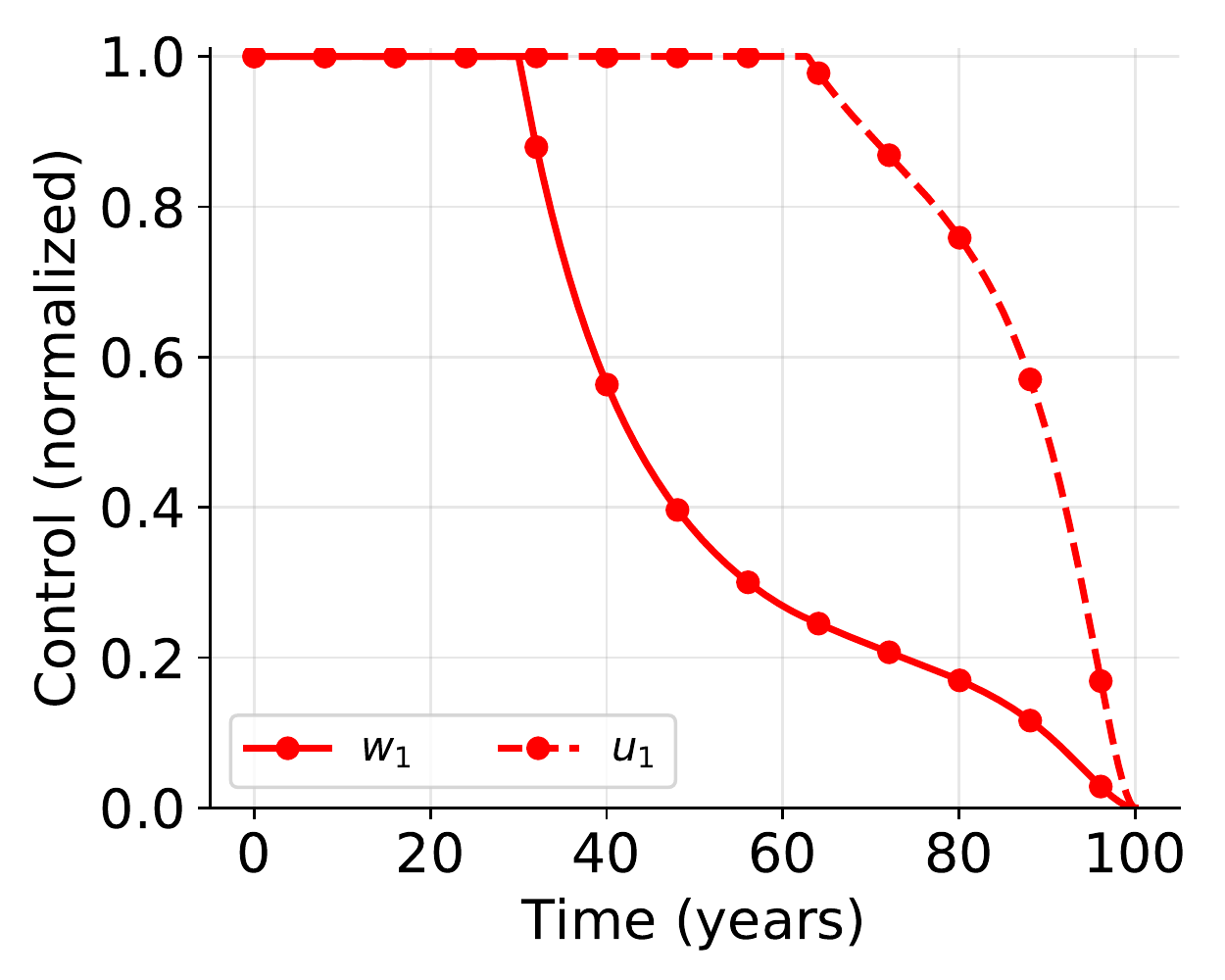}}

\subfloat[Strategy $S_{5}^{*}(t)$]{
\includegraphics[width=0.35\textwidth]{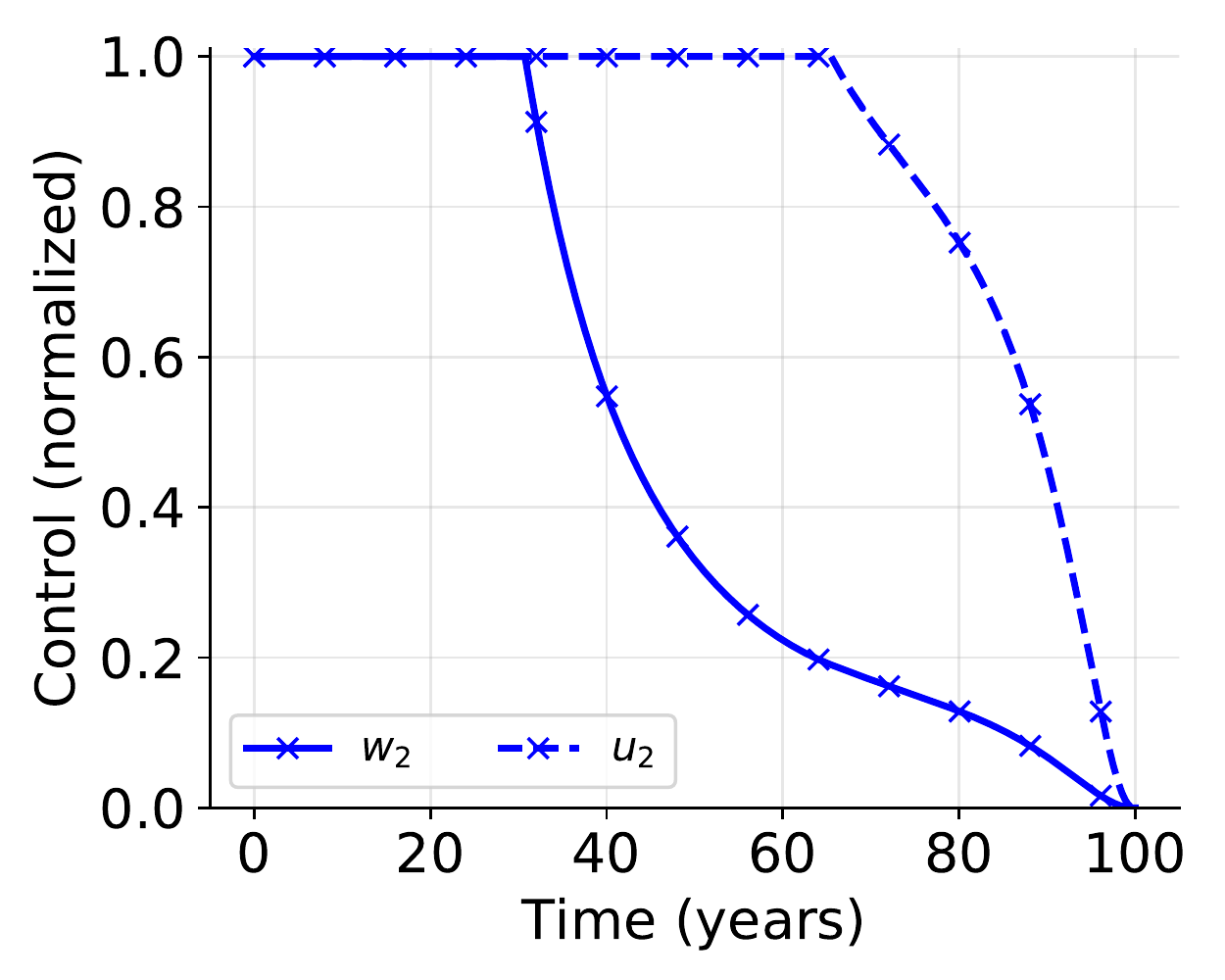}}
\subfloat[Strategy $S_{6}^{*}(t)$]{
\includegraphics[width=0.35\textwidth]{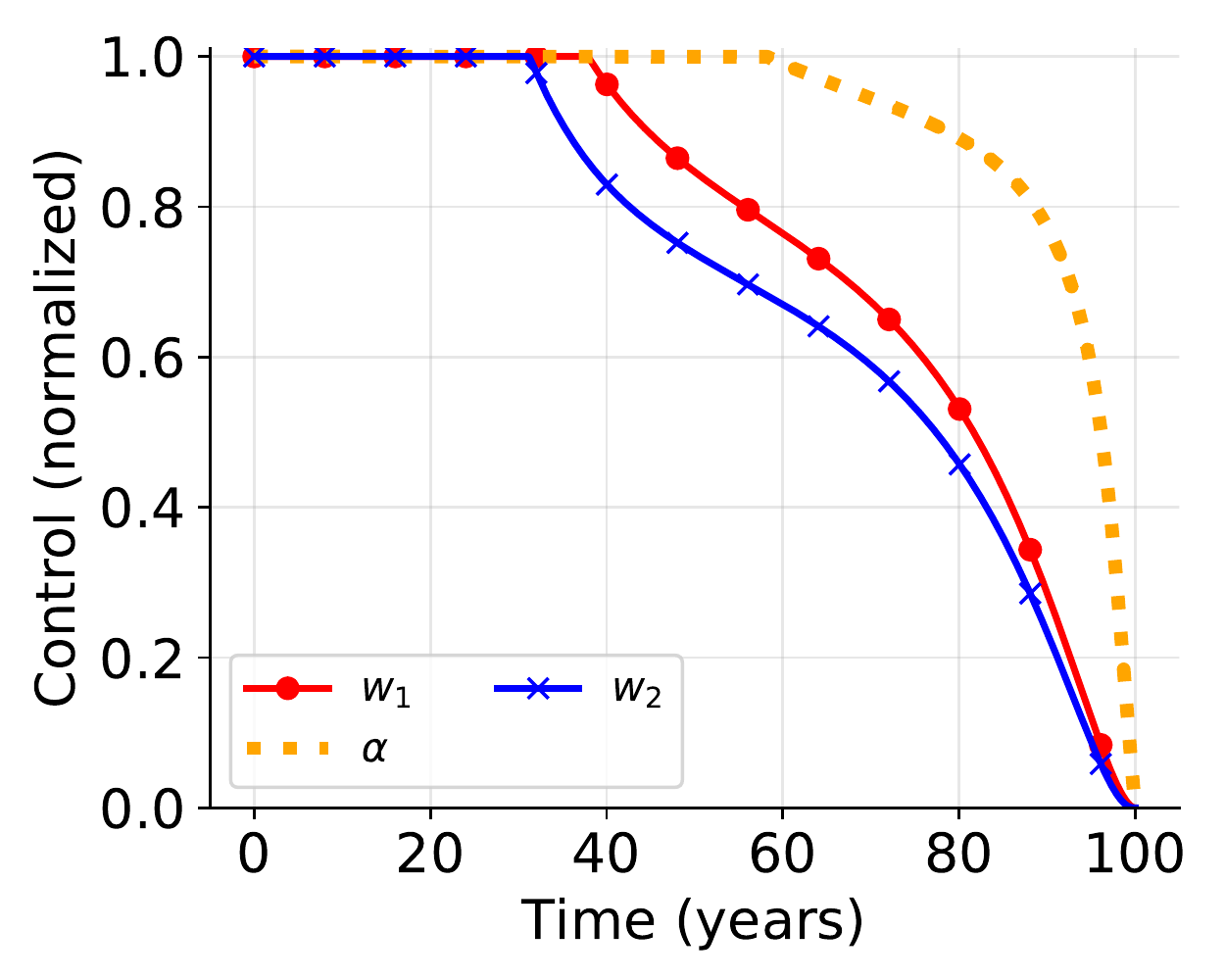}}

\subfloat[Strategy $S_{7}^{*}(t)$]{
\includegraphics[width=0.35\textwidth]{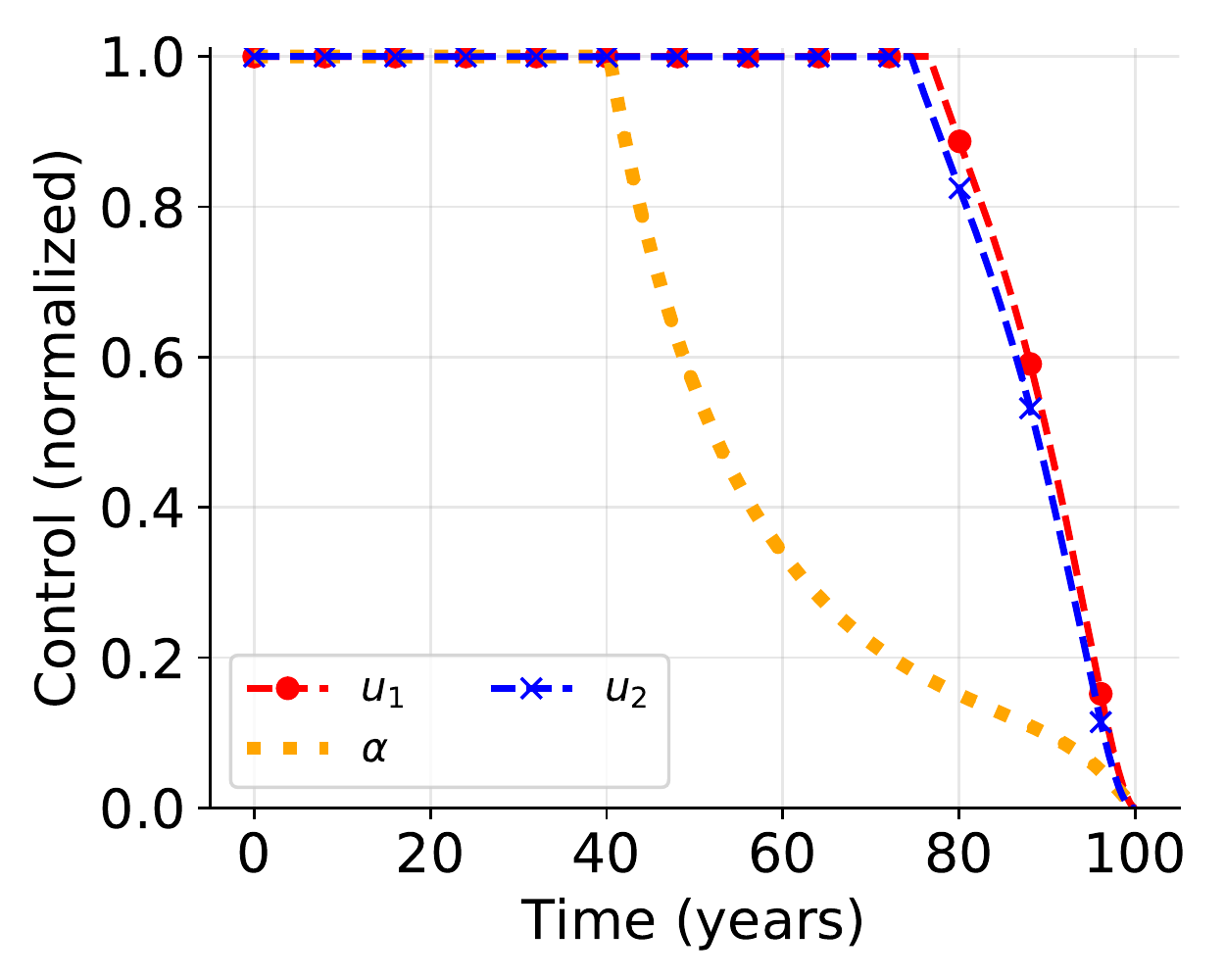}}
\subfloat[Strategy $S_{8}^{*}(t)$]{
\includegraphics[width=0.35\textwidth]{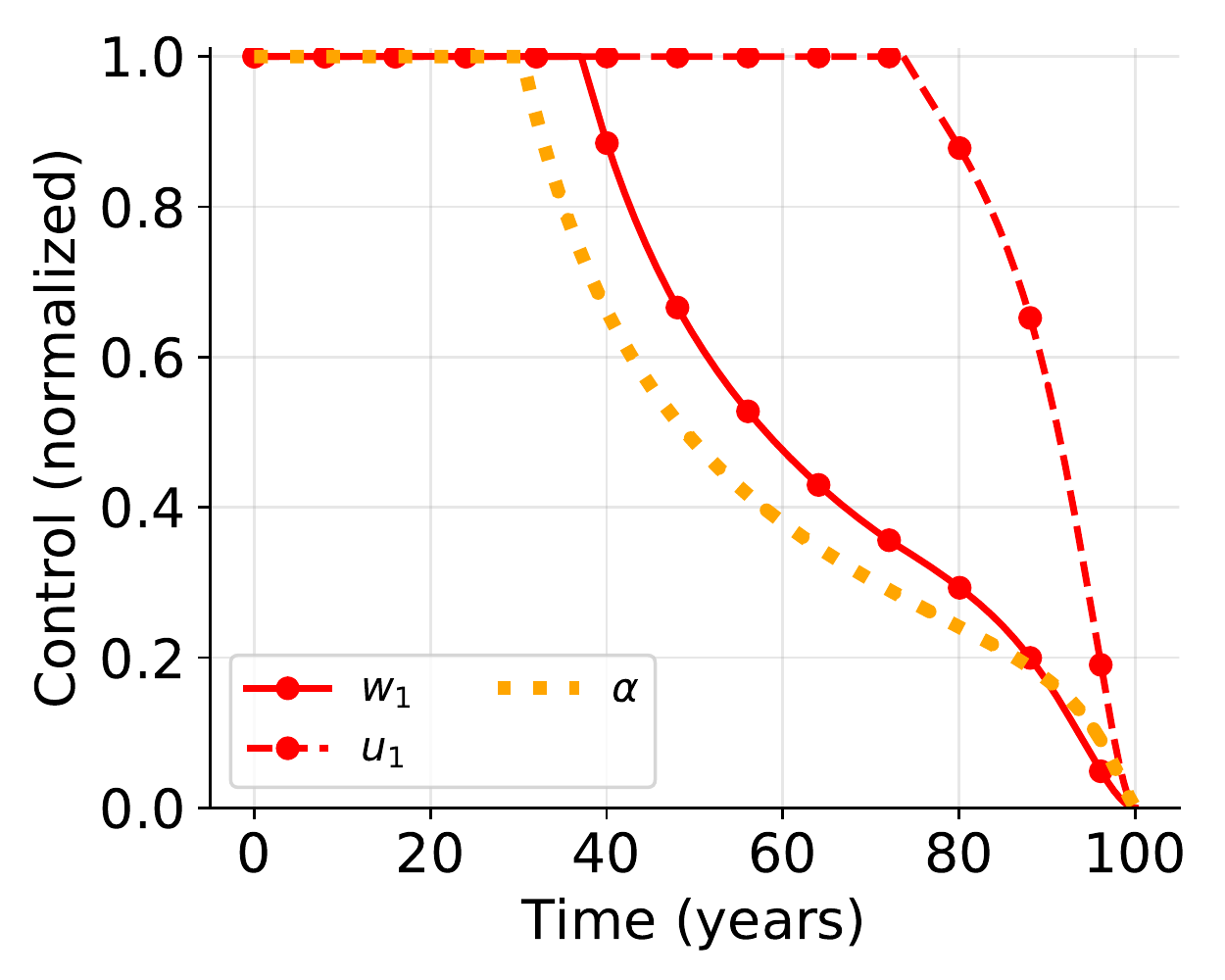}}
\caption{Time-dependent profiles of the control strategies derived from the optimal control problem \eqref{OptimalControlProblem}.}

\label{fig:optimal_profiles}
\end{figure}

The Pontryagin maximum principle converts the optimal control problem \eqref{OptimalControlProblem} into a problem of minimizing pointwise the Hamiltonian 
\begin{equation}
\begin{aligned}
H=&B_{1}U_{f}+B_{2}I_{f}+\dfrac{1}{2}\left[A_{1}(w_{1}^{2}+w_{2}^{2})+A_{2}(u_{1}^{2}+u_{2}^{2})+A_{3}\alpha_{2}\right]\\ 
&+\psi_{1}\dot{U}_{f}+\psi_{2}\dot{I}_{f}+\psi_{3}\dot{V}_{f}+\psi_{4}\dot{I}_{m}+\psi_{5}\dot{V}_{m}
\end{aligned}\label{Hamiltonian}
\end{equation}
with respect to the controls. Here, $\psi_{i}$ for $i=1,\ldots, 5$ are adjoint variables which satisfy the following system of differential equations:
\begin{equation}
\begin{aligned}
\dot{\psi}_{1}=&-B_{1}+[(1-p)\beta_{m}I_{m}+\gamma_{f}+\alpha^{*}+\mu_{f}]\psi_{1}+(p\beta_{m}I_{m}-\alpha^{*})\psi_{2}\\
&+u_{1}^{*}\psi_{3}-\beta_{f}[(1-I_{m}-V_{m})+\epsilon V_{m}]\psi_{4}+\beta_{f}\epsilon V_{m}\psi_{5},\\
\dot{\psi}_{2}=&-B_{2}+(1-p)\beta_{m}I_{m}\psi_{1}+(p\beta_{m}I_{m}+\gamma_{f}+\mu_{f})\psi_{2}+u_{1}^{*}\psi_{3}\\
&-\tilde{\beta}_{f}[(1-I_{m}-V_{m})+\epsilon V_{m}]\psi_{4}+\tilde{\beta}_{f}\epsilon V_{m}\psi_{5},\\
\dot{\psi}_{3}=&[(1-p)\psi_{1}+p\psi_{2}](1-\epsilon)\beta_{m}I_{m}+(\epsilon\beta_{m}I_{m}+u^{*}_{1}+\mu_{f}+\theta)\psi_{3},\\
\dot{\psi}_{4}=&-\beta_{m}((1-U_{f}-I_{f}-V_{f})+\epsilon V_{f})[(1-p)\psi_{1}+p\psi_{2}]+\epsilon\beta_{m}V_{f}\psi_{3}\\
&+(\beta_{f}U_{f}+\tilde{\beta}_{f}I_{f}+\gamma_{m}+\mu_{m})\psi_{4}+u^{*}_{2}\psi_{5},\\
\dot{\psi}_{5}=&(\beta_{f}U_{f}+\tilde{\beta}_{f}I_{f})((1-\epsilon)\psi_{4}+\epsilon\psi_{5})+(u_{2}^{*}+\theta+\mu_{m})\psi_{5},\\
\end{aligned}\label{adjoints}
\end{equation}
with transversality conditions $\psi_{k}(T)=0$ for $k=1,\ldots,5$. For this control problem, applying the maximum principle from, we obtain the following characterization for the optimal controls:
\begin{equation}
\begin{aligned}
w_{1}^{*}(t)&=\min\left\lbrace 1, \max\left\lbrace 0, -\dfrac{\mu_{f}}{A_{1}}\psi_{3}(t)\right\rbrace \right\rbrace ,\\
w_{2}^{*}(t)&=\min\left\lbrace 1, \max\left\lbrace 0, -\dfrac{\mu_{m}}{A_{1}}\psi_{5}(t)\right\rbrace \right\rbrace ,\\
u_{1}^{*}(t)&=\min\left\lbrace  u_{max}, \max\left\lbrace 0, -\dfrac{1}{A_{2}}(1-U_{f}(t)-I_{f}(t)-V_{f}(t))\psi_{3}(t)\right\rbrace \right\rbrace,\\
u_{2}^{*}(t)&=\min\left\lbrace  u_{max}, \max\left\lbrace 0, -\dfrac{1}{A_{2}}(1-I_{m}(t)-V_{m}(t))\psi_{5}(t)\right\rbrace \right\rbrace,\\
\alpha^{*}(t)&=\min\left\lbrace \alpha_{max}, \max\left\lbrace 0, \dfrac{1}{A_{3}}(\psi_{1}(t)-\psi_{2}(t))U_{f}(t)\right\rbrace \right\rbrace.\label{CharacterizationControls2}
\end{aligned}
\end{equation}
The control model \eqref{ControlModel}, the system of differential equations for the adjoints \eqref{adjoints} and the control characterization above form the optimality system. The optimality system will allow us to obtain time-dependent versions of the control strategies in Table \ref{tab:strategies}. To achieve this, we solve numerically such system using the forward-backward sweep method described in \citep[Chapter~4]{lenhart2007}. The time-dependent profiles for the control strategies (denoted $S_{i}^{*}(t)$, $i=1,2,\ldots,8$) are shown in Fig. \ref{fig:optimal_profiles}.

\subsection{Cost-effectiveness analysis for time-dependent control strategies}
Next, we want to determine the most cost-effective of the time-dependent control strategies $S_{i}^{*}(t)$ ($i=1,2,\ldots,8$) illustrated in Fig. \ref{fig:optimal_profiles}. As in the case of constant controls, we need to compare the costs and the effectiveness of the interventions. This is done using the cost-effectiveness algorithm presented in section \ref{sec:ICERconstant}. The results of this analysis can be found in Table \ref{tab:ICERresults_optimal}. 

\begin{table}[hbtp]
\centering
\begin{tabular}{cccc}
\toprule
Strategy &
$C(S_{i}^{*})$ &
$E(S_{i}^{*})$ &
Rank \\ 
\toprule
$S_1^{*}(t)$ 
& \$ 64.48 & 31.76 & \#8  \\ 
$S_2^{*}(t)$                
& \$ 48.65 & 30.86 & \#3 \\ 
$S_3^{*}(t)$  
& \$ 64.04 & 31.69 & \#7 \\ 
$S_4^{*}(t)$ 
& \$ 47.92 & 32.39 & \#1 \\ 
$S_5^{*}(t)$ 
& \$ 53.69 & 31.82 & \#4 \\ 
$S_6^{*}(t)$ 
& \$ 59.23 & 31.85 & \#5 \\ 
$S_7^{*}(t)$  
& \$ 64.36 & 31.99 & \#6 \\ 
$S_8^{*}(t)$ 
& \$ 50.80 & 32.61 & \#2 \\ \bottomrule
\end{tabular}
\caption{Costs (using \eqref{CostC}), cumulative level of infection averted, and rank according to the cost-effectiveness algorithm for the control strategies with time-dependent control strategies derived from \eqref{OptimalControlProblem}.}
\label{tab:ICERresults_optimal}
\end{table}
The results presented in Table \ref{tab:ICERresults_optimal} indicate that $S_{4}^{*}(t)$ is the most cost-effective strategy. This result coincides with the results for the constant control case (see Table \ref{tab:ICERresults}). Hence, female's vaccination is the intervention with the best performance for control of HPV infection. Nevertheless, for time-dependent control strategies, the second most cost-effective intervention is female's vaccination and screening, $S_{8}^{*}(t)$; whereas for constant control strategies, vaccination prior to sexual initiation, $S_{2}$, is the second most cost-effective intervention. Therefore, the rankings for the constant and the time-dependent control case differ. Moreover, it is uncertain if the time-dependent control strategies outperforms the constant control strategies as expected. Thus, we use the ICER and ACER to compare constant control strategies $S_{i}$ (see Table \ref{tab:ICERresults}) against the time-dependent control strategies $S_{i}^{*}(t)$ (see Fig. \ref{fig:optimal_profiles}).\par 
%
%
First, we compare the most cost-effective control strategies with constant and time-dependent controls, that is, $S_{4}$ versus $S_{4}^{*}(t)$. For this, we need the following computations:
\begin{equation*}
\begin{aligned}
&ACER(S_{4}^{*}(t))=\dfrac{47.92}{32.39}=1.479,\\[0.3cm]
&ICER(S_{4}^{*}(t),S_{4})=\dfrac{49.24-47.92}{32.43-32.39}=33.00.
\end{aligned}
\end{equation*}
These computations imply that $ICER(S_{4}^{*}(t),S_{4})>ACER(S_{4}^{*}(t))$. Hence, although strategy $S_4$ has higher effectiveness than strategy $S_{4}^{*}(t)$, in proportion, $S_4$ is less cost-effective than $S_{4}^{*}(t)$. In particular, the value of $ICER(S_{4}^{*}(t),S_{4})$ shows a cost saving of \$ $33$ for strategy $S_{4}^{*}(t)$ over strategy $S_{4}$. Consequently, in this case, the time-dependent control outperform the constant control.\par 
Next, we compare the second most cost-effective constant control strategy $S_2$ against the second most cost-effective time-dependent control strategy $S_8^{*}(t)$. The cost-effective ratios are:
\begin{equation*}
\begin{aligned}
&ACER(S_{2})=\dfrac{47.86}{31.04}=1.54,\\[0.3cm]
&ICER(S_{2},S_{8}^{*}(t))=\dfrac{50.80-47.86}{32.61-31.04}=1.87.
\end{aligned}
\end{equation*}
Clearly, $ICER(S_{2},S_{8}^{*}(t))>ACER(S_{2})$. Thus, the comparison shows a cost saving of \$ $1.87$ for strategy $S_{2}$ over strategy $S_{8}^{*}(t)$. Therefore, in this case, the constant control outperform the time-dependent control.\par 
Finally, we compare the third most cost-effective control strategies, that is, $S_5$ against $S_{2}^{*}(t)$. The cost-effectiveness ratios are calculated as follows: 
\begin{equation*}
\begin{aligned}
&ACER(S_{2}^{*}(t))=\dfrac{48.65}{30.86}=1.57,\\[0.3cm]
&ICER(S_{2}^{*}(t),S_{5})=\dfrac{55.07-48.65}{31.86-30.86}=6.42.
\end{aligned}
\end{equation*}
The comparison between strategies $S_5$ and $S_{2}^{*}(t)$ shows a cost saving of \$ $6.42$ for strategy $S_{2}^{*}(t)$ over strategy $S_5$. Hence, in this case, the time-dependent control outperform the constant control.\par
\begin{figure}[hbtp]
\centering
\includegraphics[width=0.5\textwidth]{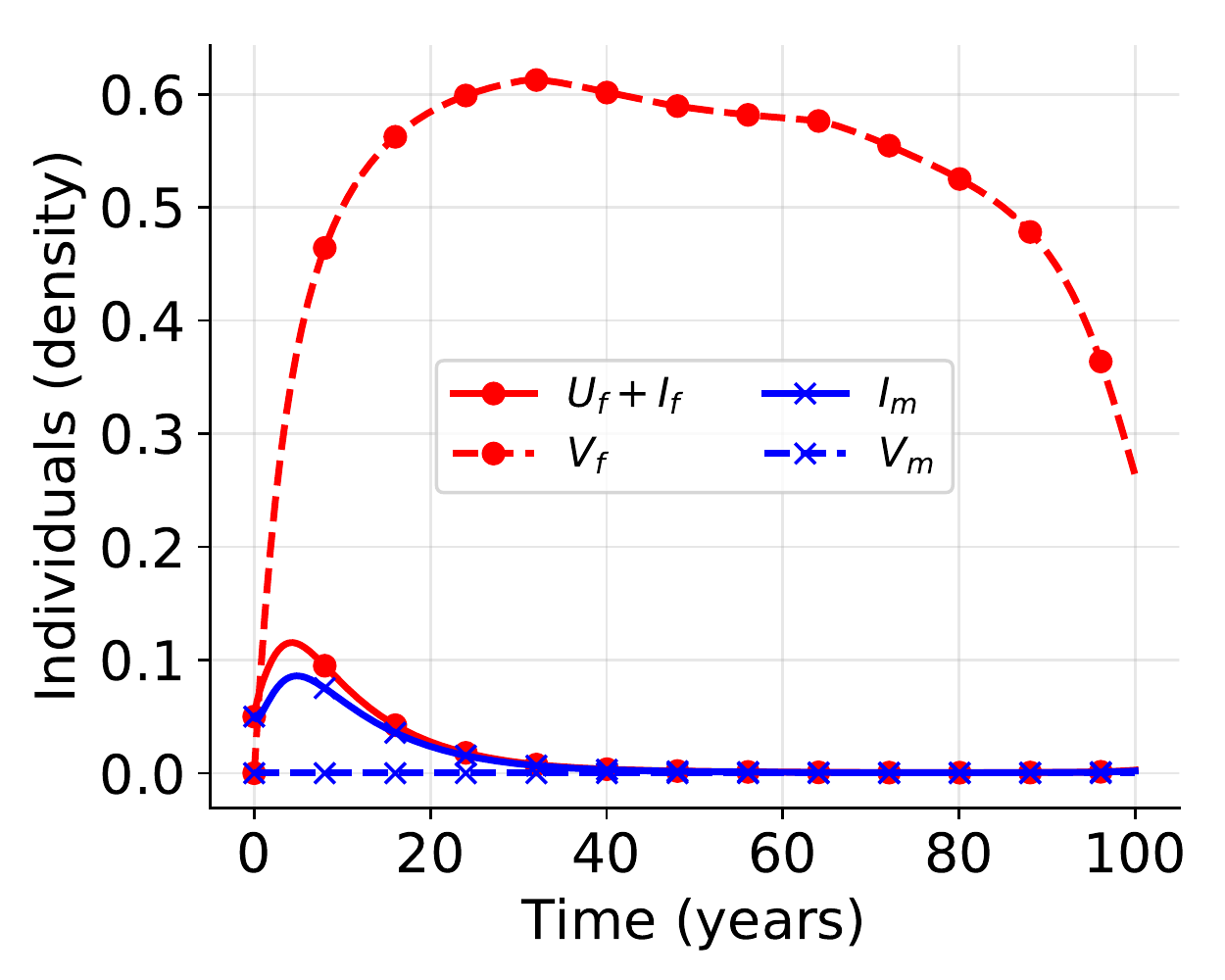}
\caption{Simulations of the HPV model showing the effects of the time-dependent control strategy $S_{4}^{*}(t)$ with optimal controls $w_1^*$ and $u_1^*$ illustrated in Fig. \ref{fig:optimal_profiles} (d).}
\label{fig:best_optimal}
\end{figure}

 
%
%
On the whole, the intervention with the best performance is the time-dependent strategy $S_{4}^{*}(t)$ (see Fig. \ref{fig:best_optimal} for the simulations of the HPV model showing the effects of strategy $S_{4}^{*}$). Nevertheless, generally speaking, time-dependent control strategies obtained by the solution of the optimal control problem are not always more cost-effective than constant control strategies. 

\section{Discussion and concluding remarks}\label{sec:discussion}
%
Data from clinical trials indicate that HPV vaccines are highly effective in preventing HPV associated diseases. Nevertheless, so far the question whether it is better, in terms of health and cost, to vaccinate females only, males only, or both sexes, remains a controversy. Some studies have shown that males' vaccination alongside females' vaccination, can be cost-effective, particularly if female vaccination coverage is low to moderate and if all potential health benefits of HPV vaccination are considered \cite{chesson2011cost}. Moreover, in circumstances where the immunization coverage in females is likely to be low, it is more optimal to vaccinate males than hard-to-reach females. However, in the majority of countries recommending HPV vaccination, the preferred strategy is vaccination of young adolescent girls aged 9-14.\par 
This work aimed to assess the cost-effectiveness of HPV healthcare interventions. We proposed constant control strategies based on combinations of the five controls incorporated in our model and used cost-effectiveness ratios to identify which strategy delivers the best effectiveness for the money invested. Comparisons of the cost-effectiveness of the proposed constant control strategies consistently show that females' vaccination, including pre-adolescent girls and adult women, is the most cost-effective strategy. Then, the next strategy with the best performance is the vaccination of school-boys and -girls before sexual initiation. Furthermore, the third most cost-effective strategy is males' vaccination. This result is unexpected because the cost function \eqref{CostC} does not depend on the presence of HPV infected males. Therefore, a very adequate choice to balance cost and protection is increasing vaccine uptake among all eligible females and extending coverage to high-risk males sub-groups.\par
Secondly, we considered the control parameters to be time-dependent and formulated an optimal control problem. We used the Pontryagin maximum principle to derive the necessary conditions for the optimal control and obtained numerical approximations via the forward-backward sweep method. This allowed us to obtain time-dependent versions of the control strategies. As in the case for constant controls, we used cost-effectiveness ratios to identify the strategy with the best performance to control HPV infection. The results confirm that the females' vaccination is the most cost-effective strategy. For this intervention, the numerical computations (see Fig. \ref{fig:optimal_profiles} (b)) indicate that initially, the vaccination rates should be applied at the maximum level and after approximately half of the time interval, these rates should gradually be reduced reaching zero at the final time.\par 
The objective functional \eqref{ObjectiveFunctional} of the optimal control problem considers the sum of weighted squares in the controls. This penalizes high levels of control administration but differs from the costs derived from function $C(S)$. Therefore, when comparing constant against time-dependent control policies, it is not known \emph{a priori} which of the strategies is the most cost-effective. Therefore, we used the ICER-ACER methodology to compare constant against time-dependent control strategies that are optimal in the sense that they minimize functional \eqref{ObjectiveFunctional}. The results indicate that time-dependent controls are not always more cost-effective than constant controls. Therefore, one must be very careful with the election of the objective functional for optimal control problems in epidemiology.


\bibliographystyle{apa}
\bibliography{bibliography_hpv}









\end{document}